\documentclass{article}

\usepackage{amsmath, amsthm, amssymb, amstext, amsfonts}
\usepackage{enumerate}
\usepackage{graphicx}
\usepackage[applemac]{inputenc}

\theoremstyle{plain}

\newtheorem{thm}{Theorem}[section]
\newtheorem{lem}[thm]{Lemma}
\newtheorem{cor}[thm]{Corollary}

\newtheorem{rem}[thm]{Remark}

\theoremstyle{definition}

\newcommand{\N}{\ensuremath{\mathbb{N}}}

\newcommand{\cT}{\ensuremath{\mathcal{T}}}

\newcommand{\sm}{\ensuremath{\setminus}}

\newcommand{\desc}[1]{\textnormal{desc}(#1)}

\newcommand{\set}[1]{\ensuremath{\{#1\}}}
\newcommand{\abs}[1]{\ensuremath{|#1|}}

\newcommand{\isom}{\ensuremath{\cong}}
\newcommand{\inv}{\ensuremath{^{-1}}}

\newcommand{\Aut}{\textnormal{Aut}}

\newcommand{\es}{\ensuremath{\emptyset}}

\newcommand{\sub}{\subseteq}

\newcommand{\comment}[1]{}

\newcommand{\nat}{{\mathbb N}}

\newcommand{\AF}{\ensuremath{\mathcal A}}
\newcommand{\BF}{\ensuremath{\mathcal B}}
\newcommand{\CF}{\ensuremath{\mathcal C}}

\newcommand{\MF}{\ensuremath{\mathcal M}}

\newcommand{\PF}{\ensuremath{\mathcal P}}

\newcommand{\SF}{\ensuremath{\mathcal S}}
\newcommand{\TF}{\ensuremath{\mathcal T}}

\begin{document}

\title{A classification of connected-homogeneous~digraphs}
\author{Matthias Hamann\and Fabian Hundertmark}
\date{Fachbereich Mathematik\\Universit\"at Hamburg\bigskip \\
\today}
\maketitle

\begin{abstract}
We classify the connected-homogeneous digraphs with more than one end. We further show that if their underlying undirected graph is not connected-homogeneous, they are highly-arc-transitive.
\end{abstract}

\section{Introduction}
A graph is called \emph{homogeneous} if every isomorphism between two finite induced subgraphs extends to an automorphism of the graph. 
If only isomorphisms between \emph{connected} induced subgraphs are required  to extend to an automorphism, the graph is called \emph{connected-homogeneous}, or simply \emph{C-homogeneous}. 
In the context of digraphs, the same notion of homogeneity and C-homogeneity applies, connectedness being taken in the underlying undirected graph. 
There are classification results for
\begin{itemize} 
\item the homogeneous graphs, \cite{C2,E,Ga1,LW, R},
\item the C-homogeneous graphs, \cite{E,Ga2, GMa, HP},
\item the homogeneous digraphs, \cite{C2,L2, L},
\end{itemize}
but not for the C-homogeneous digraphs. Our aim in this paper is to classify the C-homogeneous digraphs. Partial results towards such a classification are known for the locally finite case; they are due to Gray and Möller \cite{GMo}.

We classify the connected C-homogeneous digraphs, of any cardinality, that have more than one end.
The most important tool we use is the concept of structure trees based on vertex cut systems, introduced recently by Dunwoody and Krön \cite{DK} and used before in \cite{H,HP,K}. 
A crucial feature of this new technique is its applicability to arbitrary infinite graphs: the previously available theory of structure trees in terms of edge cuts, due to Dunwoody \cite{Du} (see also \cite{DD,Mo3,Mo2,TW}), as used by Gray and Möller \cite{GMo}, only allows for the treatment of locally finite graphs.
Our proof is based on the classification of the countable homogeneous tournaments of Lachlan \cite{L} and homogeneous bipartite graphs of Goldstern, Grossberg and Kojman \cite{GGK} and is otherwise from first principles. We reobtain the results of Gray and Möller \cite{GMo} but do not use them.
\medskip

We further study the relation between the C-homogeneous digraphs and the C-homogeneous graphs. 
A natural question arising here, is whether or not the underlying undirected graph of a C-homogeneous digraph is C-homogeneous. 
We say that a C-homogeneous digraph is of \emph{Type I} if its underlying undirected graph is C-homogeneous and otherwise it is said to be of \emph{Type~II}. Combining the results of Gray and Möller \cite{GMo} with those of Gray and Macpherson \cite{GMa} we know that there exist digraphs of both types and that there C-homogeneous graphs that do not admit a C-homogeneous orientation. 
In Section~\ref{sec_triangle} we show that connected C-homogeneous digraphs with more than one end are of Type I, if and only if they either are a tree or contain a triangle. 

Another widely studied class of digraphs are the \emph{highly-arc-transitive} digraphs, those that are $k$-arc-transitive\footnote{A (di)graph is called \emph{$k$-arc-transitive} if every (directed) path of length $k$ can be mapped to any other by an automorphism.} for all $k \in \N$. 
As a corollary of our methods, we find that the connected C-homogeneous digraphs of Type II with more than one end are highly-arc-transitive. This was previously known for locally finite such digraphs \cite{GMo}.
Unlike its undirected counterpart (cf.\ \cite{HP,TW}), the class of highly-arc-transitive digraphs is far from understood. See \cite{CPW,MMMSTZ,MMSZ,S} for papers related to highly-arc-transitive digraphs.

\section{Basics}
\subsection{Digraphs}

A digraph $D = (VD, ED)$ consists of a non-empty set $VD$, its set of \emph{vertices}, and an asymmetric (i.e.\ irreflexive and anti-symmetric) binary relation $ED$ over $VD$, its set of \emph{edges}.

We write $xy$ for an edge $(x,y) \in ED$ and say that $xy$ is directed \emph{from} $x$ \emph{to}~$y$. For $x \in VD$ we define its \emph{out-neighborhood} as $N^+(x) := \{y \in VD \:|\: xy \in ED\}$, its \emph{in-neighborhood} as $N^-(x) := \{z \in VD \:|\: zx \in ED\}$ and finally its \emph{neighborhood} as $N(x) := N^+(x) \cup N^-(x)$. Two vertices are called \emph{adjacent} if one is in the other's neighborhood. For a vertex set $X \sub VD$ the neighborhood of $X$ is defined as $N(X) := \left( \bigcup_{x \in X}N(x) \right) \sm X$ and $N^+(X), N^-(X)$ are defined analogously. For all $x\in VD$ we denote with $d^+(x),d^-(x)$ the cardinality of $N^+(x),N^-(x)$, respectively.

A sequence $x_0x_1 \dots x_k$ of pairwise distinct vertices of~$D$ with $k \in \N$ and $x_i \in N^+(x_{i -1})$ for all $1 \le i \le k$ is called a $k$-\emph{arc} from $x_0$ to $x_k$. Given two vertices $x$ and $y$ we say that $y$ is a \emph{descendant} of $x$ if there is a $k$-arc from $x$~to~$y$ for some $k \in \N$ and we define the \emph{descendant-digraph} of $x$ to be the subgraph $\desc{x} \sub D$ that is induced by the set of all its descendants. 

If $x_0x_1\dots x_n$ is a sequence of vertices such that any two subsequent vertices are adjacent then it is called a \emph{walk} and a walk of pairwise distinct vertices is called a \emph{path}. A path that is also an arc is called a {\em directed path}. A digraph is called \emph{connected} if any two vertices are joined by a path. 

A walk $x_0x_1\dots x_n$ such that $x_i \in N^+(x_{i+1}) \Leftrightarrow x_{i+1} \in N^-(x_{i+2})$ is called \emph{alternating}. If $e = xy$ and $e' = x'y'$ are contained in a common alternating walk then they are called \emph{reachable} from each other. This clearly defines an equivalence relation, the \emph{reachability relation}, on $ED$ which we denote by $\AF$, and for $e  \in ED$ we refer to the equivalence class that contains $e$ by $\AF(e)$. See also~\cite{CPW}.  

One-way infinite paths are called \emph{rays} and two rays $R_1, R_2$ are called \emph{equivalent} if for every finite vertex set $S$ both $R_1$ and $R_2$ lie eventually in the same component. This is indeed an equivalence relation on the rays of $D$ the classes of which we call the \emph{ends} of $D$. An end $\omega$ is thus a set of rays and we say that $\omega$ is \emph{contained} in a subgraph $H \subset D$ if there is a ray $R $ in $H$ such that $R \in \omega$. The same notion of an end is used for (undirected) graphs (see \cite[p.\ 202]{D}).

\subsection{Structure trees}

Let $G$ be a connected graph and let $A,B\sub VG$ be two vertex sets. The pair $(A,B)$ is a \emph{separation} of $G$ if $A\cup B= VG$ and $EG[A] \cup EG[B] = EG$.

The \emph{order} of a separation $(A,B)$ is the order of its {\em separator} $A\cap B$ and the subgraphs $G[A\sm B]$ and $G[B\sm A]$ are the \emph{wings} of $(A,B)$.
With $(A,\sim)$ we refer to the separation $(A,(VG \sm A) \cup N(VG \sm A))$.
A separation $(A,B)$ of finite order with non-empty wings is called \emph{essential} if the wing $G[A\sm B]$ is connected and no proper subset of~$A\cap B$ separates the wings of~$(A,B)$.
A {\em cut system} of~$G$ is a non-empty set $\CF$ of essential separations $(A,B)$ of~$G$ satisfying the following three conditions:
\begin{enumerate}[(i)]
\item If $(A,B)\in\CF$ then there is an $(X,Y)\in \CF$ with $X\sub B$.
\item Let $(A,B)\in\CF$ and $C$ be a component of $G[B\sm A]$. If there is a separation $(X,Y)\in\CF$ with $X\sm Y\sub C$, then the separation $(C\cup N(C),\sim)$ is also in~$\CF$.
\item If $(A,B)\in\CF$ with wings $X,Y$ and $(A',B')\in\CF$ with wings $X',Y'$ then there are components $C_1$ in $X\cap X'$ and $C_2$ in $Y\cap Y'$ or components $C_1$ in $Y\cap X'$ and $C_2$ in $X\cap Y'$ such that both $C_1$ and $C_2$ are wings of separations in~$\CF$.
\end{enumerate}

A separation $(A,B) \in \CF$ is called a $\CF$-\emph{cut}. Two $\CF$-cuts $(A_0,A_1),(B_0,B_1)$ are \emph{nested} if there are $i,j\in \set{0,1}$ such that one wing of $(A_i\cap B_j,\sim)$ does not contain any component $C$ with ${(C \cup N(C),\sim)\in\CF}$ and $A_{1-i}\cap B_{1-j}$ contains $(A_0\cap A_1)\cup(B_0\cap B_1)$.
A cut system is {\em nested} if each two of its cuts are nested.

A $\CF$-cut is {\em minimal} if there is no $\CF$-cut with smaller order. 
A {\em minimal cut system} is a cut system all whose cuts are minimal and thus have the same order.

A {\em $\CF$-separator} is a vertex set $S$ that is a separator of some separation in $\CF$.
Let $\SF$ be the set of $\CF$-separators.
A \emph{$\CF$-block} is a maximal induced subgraph $X$ of~$G$ such that
\begin{enumerate}[(i)]
\item for every $(A,B)\in\CF$ there is $VX \sub A$ or $VX\sub B$ but not both;
\item there is some $(A,B)\in\CF$ with $VX\subseteq A$ and $A\cap B\subseteq VX$.
\end{enumerate}
Let $\BF$ be the set of $\CF$-blocks.
For a nested minimal cut system $\CF$ let $\TF(\CF)$ be the graph with vertex set $\SF\cup\BF$. Two vertices $X,Y$ of~$\TF$ are adjacent if and only if either $X\in\SF$, $Y\in\BF$, and $X\sub Y$ or $X\in\BF$, $Y\in\SF$, and $Y\sub X$.
Then $\TF=\TF(\CF)$ is called the \emph{structure tree} of~$G$ and $\CF$ and by Lemma~6.2 of~\cite{DK} it is indeed a tree.

\medskip A cut system $\CF$ is called \emph{basic} if the following conditions hold:
\begin{enumerate}[(i)]
\item $\CF$ is non-empty, minimal, nested and Aut($G$)-invariant.
\item Aut($G$) acts transitively on~$\SF$.
\item For each $\CF$-cut $(A,B)$ both $A$ and $B$ contain an end of $G$ and there is no separation of smaller order that has this property.

\end{enumerate}

With the results of Dunwoody's and Krön's work on (vertex) cut systems \cite{DK} we can deduce the following theorem. 

\begin{thm}\label{thm_BasicExists}
For any graph with more than one end there is a basic cut system.\qed
\end{thm}

In the context of a digraph $D$ all concepts introduced in this section are related to the underlying undirected graph of $D$.

\subsection{Bipartite digraphs}
Let $\kappa, \lambda$ be arbitrary cardinals, and $m \in \N$. We define the \emph{directed semi-regular tree} $T_{\kappa,\lambda}$ to be the directed tree with bipartition $A \cup B$ such that $d(a) = d^+(a)~=~\kappa$ for all $a \in A$ and $d(b) = d^-(b) = \lambda$ for all $b \in B$, the \emph{complete bipartite digraph} $K_{\kappa,\lambda}$ to be the digraph with bipartition $A \cup B$ such that $|A| = \kappa$, $|B| = \lambda$ and all edges point from $A$ to $B$, the \emph{directed complement of a perfect matching} $CP_\kappa$ to be the digraph obtained from $K_{\kappa, \kappa}$ by removing a perfect matching, and the \emph{cycle} $C_{2m}$ to be the digraph obtained by orienting the undirected cycle on $2m$ vertices such that no $2$-arc arises. In the context of graphs we use the same notation to refer to the underlying undirected graph. 

We call a bipartite  graph $G$ with bipartition $X \cup Y$ \emph{generic bipartite}, if   it has the following property: For any finite disjoint subsets $U$ and $W$ of $X$ (of $Y$) there is a vertex $v$ in $Y$ (in $X$) such that $U \subseteq N(v)$ and $W \cap N(v) = \emptyset$. Any generic bipartite graph contains any countable bipartite graph as an induced subgraph, and thus up to isomorphism there is a unique countable generic bipartite graph (cp.\ \cite[p.~213]{D} and \cite[p.~98]{ES}). A \emph{generic bipartite digraph} is a digraph $D$ whose underlying undirected graph $G$ is generic bipartite with bipartition $A \cup B$ and such that all edges of $D$ are directed from $A$ to $B$.

\subsection{C-homogeneous graphs}

In order to study the C-homogeneous digraphs of Type I we make use of the classification of connected C-homogeneous graphs with more than one end from \cite{HP}, which we briefly summarize in Theorem \ref{thm_C-HomGraph}.

With $X_{\kappa,\lambda}$ we denote a graph with connectivity~$1$ such that every {\em block}, that is a maximal $2$-connected subgraph, is a complete graph on $\kappa$ vertices and every vertex lies in $\lambda$ distinct blocks.

\begin{thm}\label{thm_C-HomGraph}
A connected graph with more than one end is C-homogeneous if and only if it is isomorphic to an $X_{\kappa,\lambda}$ for cardinals $\kappa,\lambda\ge 2$.
\qed
\end{thm}

\section{Local structure}

In this chapter we summarize some preliminary results of the relation between a C-homogeneous digraph and a basic cut system $\CF$ of this digraph. In particular we investigate the local structure around $\CF$-separators.

\begin{lem}\label{lem_SeparatorEdgeless}
Let $D$ be a connected C-homogeneous digraph with more than one end.
Let $\CF$ be a basic cut system and let $S$ be a $\CF$-separator.
Then there is no edge $xy$ in~$D$ with both vertices in~$S$.
\end{lem}

\begin{proof}
Let $(A,B)\in\CF$ with $A\cap B=S$ and let us suppose that there is $xy\in ED$ with $x,y\in S$.
Let $a\in A\sm B$ and $b_0\in B\sm A$.
Let us first consider the case that $ay,b_0y\in ED$.
Then there is an automorphism $\alpha_0$ of~$D$ such that $(xy)^{\alpha_0}=b_0y$.
Hence there is another vertex $b_1\in B$ with $b_1y\in ED$ and such that $b_1$ and $a$ are separated from each other by $S$ and $S^{\alpha_0}$.
By repeating this process with an $\alpha_i$ that maps $b_{i-1}$ onto $b_i$ we get a further vertex $b_{i+1}$ that is separated by $(A\cap B)^{\alpha_0\ldots\alpha_i}$ from~$a$.
After the step $\abs{A\cap B}$ there has to be some $b_i$ that is also separated from $b:=b_{\abs{A\cap B}+1}$ by $(A\cap B)^{\alpha_0\ldots\alpha_{\abs{A\cap B}}}$.
But then each $\CF$-separator that separates $b_i$ from $b$ also has to separate $a$ from $b$.
Since $S$ separates $a$ and $b$ but not $b_i$ and $b$ and since $D[a,y,b]$ is isomorphic to $D[b_i,y,b]$, we get a contradiction to the C-homogeneity of~$D$.

So let us assume that there are vertices $a\in A\sm B$, $b\in B\sm A$ with $by,ya\in ED$.
Let $\alpha$ be an automorphism of~$D$ with $(xy)^\alpha=ya$.
Then there is a neighbor $c$ of $y$ that is separated from $b$ by $(A\cap B)^\alpha$.
If $cy\in ED$ then we may take the vertices $c,b$ instead of $a,b$ and get a contradiction by the first case above.
Thus we suppose that $yc\in ED$.
But then we can map the digraph $D[b,y,a]$ onto $D[b,y,c]$.
Since every separator that separates the vertices $b$ and $a$ also separates $b$ and $c$ but $(A\cap B)^\beta$ separates $b$ and $c$ but not $b$ and $a$, we get a contradiction.
The case $ay,yb\in ED$ is analog.

Let us finally assume that there are vertices $a\in A\sm B$ and $b\in B\sm A$ such that $ya,yb\in ED$.
By considering the digraph $D\inv$ instead of~$D$ we also may assume that there are $a'\in A\sm B$ and $b'\in B\sm A$ with $a'x,b'x\in ED$.
Let $\alpha$ be an automorphism of~$D$ with $(xy)^\alpha=yb$.
Then there is a vertex $b'\in B\sm A$ that is separated by $(A\cap B)^\alpha$ from $a$ and such that $b'y\in ED$.
But then we have the situation of the previous case and thus we know that no such edge $xy$ exists.
\end{proof}

\begin{lem}\label{lem_No2Arc}
Let $D$ be a connected C-homogeneous digraph with more than one end and let $\CF$ be a basic cut system. Then for each $2$-arc $P$ in $D$ we have $\left|P \cap S \right| \le 1$ for all $\CF$-separators $S$.
\end{lem}
\begin{proof}
Let $P = xay$ be a $2$-arc in $D$ and $S$ a $\CF$-separator. By Lemma~\ref{lem_SeparatorEdgeless} we only have to show that $S$ cannot contain both $x$ and $y$. So assume $\{x,y\} \subseteq S$. Let $(A,B)\in\CF$ with $A\cap B=S$ and $a \in A$. Since $D$ is transitive there is an arc $zx$ in $D$. If $z$ lies in $A$ consider a neighbor $z^\prime$ of $x$ in $B$. Now either $zxa$, $zxz^\prime$ or $z^\prime xa$ is an induced $2$-arc in $D$, which we denote by $Q$, with one vertex in $A \sm B$ and one vertex in $B \sm A$. Because $D$ is connected-homogeneous there is an automorphism $\alpha$ with $P^\alpha = Q$. But then $S^\alpha$ contains  vertices of both wings of $(A,B)$, contradicting the nestedness of $\CF$.
\end{proof}

\begin{lem}\label{lem_NokArc}
Let $D$ be a connected C-homogeneous digraph with more than one end, let $\CF$ be a basic cut system of~$D$, and let $S$ be a $\CF$-separator.
Then there is no directed path in~$D$ with both endvertices in~$S$.
\end{lem}

\begin{proof}
Suppose that there is such a path $P$.
We may choose the path such that it has minimal length.
Then all of the vertices of~$P$ lie in the same $\CF$-block $X$.
By Lemma~\ref{lem_No2Arc} the endvertices of any directed path of length $2$ are separated by a $\CF$-separator.
Hence no directed path of length at least $2$ can lie in any $\CF$-block.
\end{proof}

\begin{lem}\label{lem_NoInAndOut}
Let $D$ be a connected C-homogeneous triangle-free digraph with more than one end, and let $\CF$ be a basic cut system. Then for any cut $(A,B) \in \CF$ there is no path $xyz$ in $D\left[A\right]$ with $y \in A\cap B$.
\end{lem}
\begin{proof} By Lemma \ref{lem_SeparatorEdgeless} we only have to show that given a cut $(A,B) \in \CF$ there is no path $xyz$ in $D$ such that $y \in S := A \cap B$ and $x,z \in A \sm B$. So let us suppose there is such a path.
Then $y$ has a neighbor $b \in B \sm A$. We may assume that their connecting edge is pointing towards $y$, since otherwise changing the direction of each edge gives a digraph $D^\prime$ which is C-homogeneous and has this property.

Suppose that there is a second neighbor $c\in B\sm A$ of~$y$. If there is $yc\in ED$ then there is an $\alpha\in \Aut(D)$ that fixes $b,y,z$ and with $x^\alpha=c$, $c^\alpha=x$. But then the separations $(A,B)$ and $(A^\alpha,B^\alpha)$ are not nested. Thus we may assume that $cy\in ED$. In this situation let $\beta$ be an automorphism of~$D$ that fixes $x,y,b$ and maps $z$ onto $c$ and vice versa---a contradiction as before.

So $b$ is the unique neighbor of $y$ in $B$. We may assume that there is another vertex $a$, say, that lies in $S$, since otherwise $y$ would seperate $x$ from $z$, contradicting the fact that $x$ and $z$ lie in the same component of $D - S$. Now consider a path $P$ in $D$ connecting $a$ and $y$ and let $\cT$ denote the structure tree of~$D$ and~$\CF$. Let $\MF$ be the set of $\CF$-blocks containing edges of $P$. 
Since $\CF$-separators do not contain any edge, distinct blocks cannot contain a common edge.
Thus we choose a block $M \in \MF$ whose distance to~$S$ in $\cT$ is maximal with respect to~$\MF$.

Now each nontrivial component of $P \cap M$ has to contain exactly two edges: An isolated edge would either be contained in a separator, in contradiction to Lemma \ref{lem_SeparatorEdgeless}, or it would connect $M$ to two distinct neighbors in $\cT \cap \MF$, contradicting the choice of $M$. If there is a segment of $P$ in $M$ with a length of at least three, then it contains either a directed subsegment, isomorphic to $byz$, or a subsegment isomorphic to $by \cup xy$. In each case there exists an isomorphism $\varphi$ such that $S^\varphi$ separates the endvertices of this subsegment, which is impossible since $M$ is a $\CF$-block.

Considering an arbitrary nontrivial component of $P \cap M$, its two edges have a common vertex which we denote by $m$. With an analog argument as above, both edges are directed away from $m$. Let us denote their head by $u$ and $v$, respectively. By construction, $u$ and $v$ lie both in the separator $S_M \subset M$ that lies on the unique shortest path between $M$ and $S$ in~$\cT$. Consider an arbitrary cut with seperator $S_M$. Then $u$ has a neighbor $u^\prime$ in the wing not containing $m$. Let $\psi$ be an automorphism with $(mu)^\psi=by$ and either $(uu^\prime)^\psi=yz$, if $uu^\prime \in ED$ or $(u^\prime u)^\psi=xy$, if $u^\prime u \in ED$. Since $\CF$ is nested we have $S_M^\psi \subset B$ which means that $x$ and $z$ are seperated from $b$ by $S_M^\psi$. By relabeling $S := S_M^\psi$ and $a := v^\psi$, if neccessary, we may assume that $b$ sends an edge to $a$. 

Then there is a neighbor $z^\prime$ of $a$ in $A \sm B, $ and we can find an automorphism $\gamma$ with $(by)^\gamma = ba$ and either $x^\gamma = z^\prime$ or $z^\gamma = z^\prime$, depending on the orientation of the edge between $b$ and $z^\prime$. Again by the nestedness of $\CF$ we have $S^\gamma \subset B$ and also $B^\gamma \sub B$. And since $x$ is seperated from $b$ by $S^\gamma$ we have $y \in S^\gamma$. But that implies that $y$ and $a$ both have $b$ as their unique neighbor in $B^\gamma$. Hence, $S^\gamma \sm \{y,a\} \cup \{b\}$ is a seperator in $D$ that seperates ends and has smaller cardinality, contradicting the fact that $\CF$ is basic.
\end{proof}

\begin{lem}\label{lem_OneOut/InBlock}
Let $D$ be a connected C-homogeneous triangle-free digraph that is not a tree and that has more than one end, and let $\CF$ be a basic cut system of~$D$.
Let $S$ be a $\CF$-separator and let $s\in S$.
Then there is precisely one $\CF$-block that contains $s$ and all edges directed away from $s$, and there is precisely one $\CF$-block that contains $s$ and all edges directed towards $s$.
Furthermore there is $d^+(s)>1$ and $d^-(s)>1$.
\end{lem}

\begin{proof}
By Lemma~\ref{lem_NoInAndOut} there is at most one kind of neighbors in each $\CF$-block.
Suppose first that there is one $\CF$-block with only one neighbor $a$ of~$s$.
We may assume that $as\in ED$.
By C-homogeneity each $\CF$-block $Y$ that contains an in-neighbor of~$s$ contains no other neighbor of~$s$.
Thus each component of each $\CF$-block is either a single vertex or a star the edges of which are directed towards the center of the star.
Then, since $D$ is not a tree, there is a second vertex $t\in S$.
In every component $C$ of $D-S$ there is an (undirected) $s$-$t$-path $P$.
Let $X$ be a $\CF$-block with maximal distance to~$S$ in~$\TF$ such that there is at least one edge from $P$ in~$X$.
By Lemma~\ref{lem_SeparatorEdgeless} there is at least a second edge in~$X$.
As each component of~$X$ that contains edges is a star, the longest subpath of~$P$ that lies completely in~$X$ has length $2$.
Let $xyz$ be such a subpath.
Then $xy,zy\in ED$ and $y$ is the only neighbor of each $x$ and $z$ in~$X$.
Let $S'$ be that $\CF$-separator that contains $x$ and $z$ with $S'\sub X$.
By replacing $S'$ with $S'\sm\{x,z\}\cup\{y\}$ we get a contradiction as in the proof of the previous lemma.

Thus a $\CF$-block cannot contain $s$ together with a single neighbor of~$s$ and by C-homogeneity there has to be one $\CF$-block that contains all in-neighbors of~$s$ and one that contains all out-neighbors of~$s$.
\end{proof}

\begin{lem}\label{lem_SeparatorDegree2}
Let $D$ be a connected C-homogeneous triangle-free digraph that is not a tree and that has more than one end, and let $\CF$ be a basic cut system of~$D$.
Then each $\CF$-separator has degree two in the structure tree $\cT$ for~$D$ and~$\CF$.
\end{lem}
\begin{proof}
Let $S$ be a $C$-separator. Then for each component $X$ of $\cT - S$ the vertex set $(\bigcup X)\sm S$ is the union of components of $D - S$. Since each $s \in S$ has a neighbor in each component of $D - S$, it also has at least one neighbor in each component of $\cT - S$. With Lemma \ref{lem_OneOut/InBlock} we have $d_{\cT}(S) = 2$.
\end{proof}

\noindent If we combine Lemma \ref{lem_OneOut/InBlock} and Lemma \ref{lem_SeparatorDegree2} we get the following
\begin{cor}\label{cor_SeparatorToBlock}
Let $D$ be a connected C-homogeneous triangle-free digraph that is not a tree and that has more than one end, and let $\CF$ be a basic cut system of~$D$. Let $B$ be a $\CF$-block, $S \subset B$ a $\CF$-separator and $s \in S$. If $s$ has no neighbor in $B$, then there is exactly one $\CF$-separator $S^\prime \subset B$ such that $s \in S^\prime \cap S$. If $s$ has a neighbor in $B$, then $S$ is the only $\CF$-separator in $B$ that contains $s$. \qed
\end{cor}

\begin{lem}\label{lem_WithTrianglesBlocksComplete}
Let $D$ be a connected C-homogeneous digraph with more than one end that embeds a triangle, and let $\CF$ be a basic cut system of~$D$. Then every $\CF$-block that contains edges is complete.
\end{lem}

\begin{proof}
Let $S$ be a $\CF$-separator and let $x\in S$.
Then $x$ has adjacent vertices in both wings of each cut $(A,B)\in\CF$ with $A\cap B=S$.
As $D$ contains triangles, each edge lies on a triangle.
We know that each wing of $(A,B)$ contains both an in- and an out-neighbor of~$x$, as any triangle contains a $2$-arc and $D$ is edge-transitive.
Thus every induced path of length $2$ in~$D$ can be mapped on a path crossing $S$, i.e.\ a path both end vertices of which lie in distinct wings of $(A,B)$.
Hence no two vertices in the same $\CF$-block can have distance $2$ from each other and, in particular, every component of every $\CF$-block has diameter~$1$.

To prove that each $\CF$-block has diameter~$1$ we just have to show that each $\CF$-block is connected.
So let us suppose that this is not the case.
Let $X$ be a $\CF$-block and let $P$ be a minimal (undirected) path in~$D$ from one component of~$X$ to another.
Let $Y$ be a $\CF$-block with maximal distance in the structure tree of~$D$ and $\CF$ to~$X$ that contains edges.
By Lemma~\ref{lem_SeparatorEdgeless} the block $Y$ has to contain at least two edges and there are two non-adjacent vertices in the same component of~$Y$.
This contradicts that these components are complete graphs.
Hence each $\CF$-block that contains edges has precisely one component which has diameter~$1$.
\end{proof}

\section{C-homogeneous digraphs of Type I}\label{sec_triangle}

In this section we will completely classify the countable connected C-homoge\-neous digraphs of Type I with more than one end and give - apart from the classification of infinite uncountable homogeneous tournaments - a classification of uncountable such digraphs.
As a part of the countable classification we apply a theorem of Lachlan \cite{L}, see also \cite{C}, on countable homogeneous tournaments.
Lachlan proved that there are precisely $5$ such tournaments.
Three of them are infinite, one is the digraph on one vertex with no edge and one is the directed triangle.
For the uncountable case there is up to now no such classification of homogeneous tournaments.

To state Lachlan's theorem let us first define the countable tournament $\PF$ to be the digraph with the rationals in the intervall $[-\pi,\pi]$ as vertex set and direct the edge from $x$ to~$y$ if $$x-y\le\pi\mod 2\pi$$ and from $y$ to~$x$ otherwise.
The generic countable tournament is the unique (cp.\ \cite[p.~213]{D}, and \cite[p.~98]{ES}) countable homogeneous tournament that embeds all finite tournaments.

\begin{thm}[{\cite[Theorem~3.6]{C}}]\label{thm_Lachlan}
There are up to isomorphism only $5$ countable homogeneous tournaments: the trivial tournament on one vertex, the directed triangle, the generic tournament on $\omega$ vertices, the tournament that is isomorphic to the rationals with the usual order, and the tournament $\PF$ described above.\qed
\end{thm}

For a homogeneous tournament $T$ let $X_\lambda(T)$ denote the digraph that has connectivity $1$ and each block is isomorphic to~$T$ and each vertex is a cut vertex and lies in $\lambda$ distinct copies of~$T$.
Thus the underlying undirected graph is a distance-transitive graph as described in~\cite{HP,Ma,Mo}.

\begin{thm}\label{thm_triangles}
Let $D$ be a connected digraph with more than one end.
Then $D$ is C-homogeneous of Type I if and only if  one of the following statements holds:
\begin{enumerate}[(1)]
\item $D$ is a tree with constant in- and out-degree;
\item $D$ is isomorphic to a $X_\lambda(T^\kappa)$, where $\kappa$ and $\lambda$ are cardinals with $\lambda\ge 2$ and $\kappa$ either $3$ or infinite and $T^\kappa$ is a homogeneous tournament on $\kappa$ vertices.
\end{enumerate}
\end{thm}

\begin{proof}
Let us first assume that $D$ is a C-homogeneous digraph of Type I. Then the underlying undirected graph is isomorphic to a $X_{\kappa,\lambda}$ for cardinals $\kappa,\lambda\ge 2$.
If $\kappa=2$, then $D$ is a tree with constant in- and out-degree, so we may assume $\kappa\ge 3$.
As each block has to be homogeneous, we conclude from Theorem~\ref{thm_Lachlan} that the cardinal $\kappa$ has to be either $3$ or infinite.
This proves the necessity-part of the statement.

As the digraphs of part (1) are obviously C-homogeneous of Type I, we just have to assume for the remaining part that $D$ is isomorphic to $X_\lambda(T^\kappa)$ for a cardinal $\lambda\ge 2$ and a homogeneous tournament $T^\kappa$ on $\kappa$ vertices for a cardinal $\kappa$ that is either $3$ or infinite.
Let $\CF$ be a basic cut system of~$D$.
Let $X$ and $Y$ be two connected induced finite and isomorphic subdigraphs of~$D$.
Let $\varphi$ be the isomorphism from $X$ to~$Y$.
If $X$ has no cut vertex, then $X$ lies in a subgraph of~$D$ that is a homogeneous tournament and the same is true for $Y$, so $\varphi$ extends to an automorphism of~$D$.
So let $x\in VX$ be a cut vertex of~$X$. Hence $x^\varphi$ is a cut vertex of~$Y$.
It is straight forward to see that for any $\CF$-block $B$ the image of $X\cap B$ in~$Y$ is precisely the intersection of $Y$ with a $\CF$-block $A$.
Since the $\CF$-blocks are all isomorphic homogeneous tournaments, the isomorphism from $X\cap B$ to $Y\cap A$ extends to an isomorphism from $X$ to~$Y$.
Thus the isomorphism from $X$ to~$Y$ easily extends to an automorphism of~$D$.
Since the underlying undirected graph is C-homogeneous by Theorem~\ref{thm_C-HomGraph}, $D$ is C-homogeneous of Type I.
\end{proof}

Lachlan's theorem together with Theorem~\ref{thm_triangles} enables us to give a complete classification of countable connected C-homogeneous digraphs of Type I and with more than one end:

\begin{cor}\label{cor_Triangle}
Let $D$ be a countable connected digraph with more than one end.
Then $D$ is C-homogeneous of Type I if and only if one of the following assertions holds:
\begin{enumerate}[(1)]
\item $D$ is a tree with constant countable in- and out-degree;
\item $D$ is isomorphic to a $X_\lambda(Y)$, where $\kappa$ is a countable cardinal greater or equal to~$2$ and $Y$ is one of the four non-trivial homogeneous tournaments of Theorem~$\ref{thm_Lachlan}$.\qed
\end{enumerate}
\end{cor}

\section{Reachability and descendant digraphs}

In this section we prove that, if a connected C-homogeneous digraph $D$ with more than one end contains no triangles, then $D$ is highly-arc-transitive, each reachability digraph of $D$ is bipartite, and, if furthermore $D$ has infinitely many ends, then the descendants of each vertex in $D$ induce a tree.
All these properties were proved to be true in the case that $D$ is locally finite, see~\cite[Theorem~4.1]{GMo}.

\begin{thm}\label{thm_HighlyArcTransitive}
Let $D$ be a connected C-homogeneous triangle-free digraph with more than one end. Then $D$ is highly-arc-transitive.
\end{thm}
\begin{proof}Let $\CF$ be a basic cut system. It suffices to show that each directed path is induced. Suppose this is not the case. Then there is a smallest $k$ such that there is a $k$-arc $A = x_0 \dots x_k$ that is not induced. Hence there is an edge between $x_0$ and $x_k$. Consider a $\CF$-separator $S$ that contains $x_1$. By Lemma \ref{lem_NokArc} we have $x_k \notin S$; hence $x_0$ and $x_k$ lie on the same side of $S$. But then the same holds for $x_{k-1}$ and so on. So finally $x_0$ and $x_2$ have to lie on the same side of $S$, in contradiction to Lemma \ref{lem_NoInAndOut}.
\end{proof}

In an edge-transitive digraph all \emph{reachability digraphs} $\Delta_e := D[\AF(e)]$ with $e\in ED$ are isomorphic, so we may denote a representative of their isomorphism type by $\Delta(D)$.
Furthermore Cameron, Praeger and Wormald \cite[Proposition~1.1]{CPW} proved that the reachability relation in such a digraph is either universal or the corresponding reachability digraph is bipartite.
We will now prove that the reachability relation is not universal in our case.

\begin{thm}\label{thm_DescendentReachability}
Let $D$ be a connected C-homogeneous triangle-free digraph with more than one end.
Then $\Delta(D)$ is bipartite and if D is not a tree, then each $\Delta_e$ with $e\in ED$ is a component of a $\CF$-block. Furthermore, if $D$ has infinitely many ends, then every descendant digraph $\desc{x}$ with $x \in VD$ is a tree.
\end{thm}

\begin{proof}
Let $\CF$ be a basic cut system.
We first show that either $D$ is a tree or any $\Delta_e$ with $e\in ED$ is a component of a $\CF$-block.
Let us assume that $D$ is not a tree.
If a vertex $x$ has at most one neighbor in a $\CF$-block~$X$ with $x\in X$, then the digraph is a tree by Lemma~\ref{lem_OneOut/InBlock}.
Thus no separator can separate two vertices $x,y$ for which there is $z\in VD$ with $xz,yz\in ED$ or with $zx,zy\in ED$.
Thus each $\Delta_e$ lies in a $\CF$-block.
As there are induced paths of length $2$ crossing some $\CF$-separator and as $D$ contains no triangle, a component of a $\CF$-block $X$ cannot contain more vertices than $\Delta_e$ with $e\in E(D[X])$ contains. Thus $\Delta_e$ is a component of a $\CF$-block.

Suppose that $\Delta(D)$ is not bipartite.
Then there is a cycle of odd length in $\Delta(D)$.
Thus there has to be a directed path of length at least $2$ on that cycle.
By Lemma~\ref{lem_No2Arc} this path lies in distinct $\CF$-blocks.
This is not possible as shown above and thus $\Delta(D)$ has to be bipartite.

Now suppose that there is $x \in VD$ such that $\desc{x}$ contains a cycle. So by transitivity there is a descendant $y$ of $x$ such that there are two $x$-$y$-arcs that are apart from $x$ and $y$ totally disjoint. Thus, since we are C-homogeneous, any two out-neighbors of $x$ have a common descendant. Assume that there are two distinct $\CF$-separators $S, S^\prime$ such that both $Y := S \sm S^\prime$ and $Y^\prime := S^\prime \sm S$ contain an out neighbor of~$x$. Then it exists a vertex $z$ in $D$ with $Y$-$z$- and $Y'$-$z$-arcs. But by the Lemmas~\ref{lem_NokArc} and \ref{lem_NoInAndOut} the vertices $x$ and $z$ cannot lie on the same side of $S$ and $S^\prime$, respectively, hence $S$ and $S^\prime$ meet on both sides, a contradiction to the nestedness of $\CF$. Thus there is a $\CF$-separator $S_{+1}$ that contains the whole out-neighborhood of $x$. This implies that all descendants of distance $k$ are contained in a common $\CF$-separator $S_{+k}$, since either all distinct $k$-arcs originated at $x$ are disjoint, and we can apply the same argument as above, or each two of those $k$-arcs intersect in a vertex $x^\prime$ in $D$ that has the same distance to $x$ on both arcs by Lemma~\ref{lem_NokArc}, and we are home by induction. 

With a symmetric argument we get that  each $k$-arc that ends in $x$ has to start in a common $\CF$-separator $S_{-k}$. For a path $P$ in $D$ that starts in $x$, let $\sigma(P)$ denote the difference of the number of edges in $P$ that are directed away from $x$ (with respect to $P$) minus the number of edges of the other type. Then one easily checks that the endvertex of $P$ lies in $S_{\sigma(P)}$. Since all $\CF$-separators have the same finite order~$s$, say, there can be at most $2 s$ rays that are eventually pairwise disjoint. Hence $D$ has finitely many ends, which proves the last statement of the theorem.
\end{proof}

\begin{lem}\label{lem_SeparatorAndReachability}
Let $D$ be a connected C-homogeneous triangle-free digraph with more than one end and let $\CF$ be a basic cut system of~$D$.
Then for each $\CF$-separator $S$ of order at least $2$ there is a reachability digraph $\Delta_e$ and a $\CF$-block $K$ such that $|S\cap \Delta_e|\ge 2$, $\Delta_e\sub K$, and $S\sub K$.
\end{lem}

\begin{proof}
Let $S$ be a $\CF$-separator with $|S| \ge 2$.
Suppose that there is no reachability digraph $\Delta_e$ with $|S\cap \Delta_e|\ge 2$.
Let $x,y\in S$ and let $P$ be an $x$-$y$-path in a component of $D-S$.
Let $B$ be a $\CF$-block that contains edges of~$P$ and such that $d_\TF(S,B)$ is maximal with this property.
Then the $\CF$-separator $S_B \sub B$ that separates $S$ and $B$ in~$\TF$ has the desired property and thus each $\CF$-separator has it, in contradiction to the assumption.
\end{proof}

We have roughly described the global structure of C-homogeneous digraphs. To investigate the local structure of these graphs, we show that the underlying undirected graph of each reachability digraph is a connected C-homogeneous bipartite graph. Such graphs shall be described in the next section.

\begin{lem}\label{lem_ReachabilityCHomBipartite}
Let $D$ be a triangle-free connected C-homogeneous digraph with more than one end. Then the underlying undirected graph of $\Delta(D)$ is a connected C-homogeneous bipartite graph.
\end{lem}
\begin{proof}
By Theorem \ref{thm_DescendentReachability} $\Delta(D)$ is bipartite. The remainder of the proof is the same as the proof of the locally finite case in \cite[Lemma 4.3]{GMo}. 
\end{proof}

\section{C-homogeneous bipartite graphs}\label{sec_bipartite}
In this chapter we complete the classification of connected C-homogeneous bipartite graphs, which was already done for locally finite graphs, by Gray and Möller \cite{GMo}. They already mentioned that their work should be extendable with not too much effort -- and indeed this section has essentially the same structure. 

The proof of the locally finite analog \cite[Lemma 4.4]{GMo} of Lemma~\ref{lem_EmbeddingC4} is self contained and does not use the local finiteness of the graph. Thus we can omit the proof here.

\begin{lem}\label{lem_EmbeddingC4}
Let $G$ be a connected C-homogeneous bipartite graph with bipartition $X \cup Y$. If $G$ is not a tree and has at least one vertex with degree greater than $2$ then $G$ embeds $C_4$ as an induced subgraph.\qed
\end{lem} 

Let $G$ be a bipartite graph with bipartition $X \cup Y$. Then for each edge $\{x,y\} \in EG$ we define the neighborhood graph to be:
$$\Omega(x,y) := G[N(x) + N(y) - \{x,y\}]$$
A C-homogeneous graph $G$ is, in particular, edge-transitive, hence there is a unique neighborhood graph $\Omega(G)$.

\begin{lem}\label{lem_Omega}
Let $G$ be a connected C-homogeneous bipartite graph. Then $\Omega(G)$ is a homogeneous bipartite graph, and therefore is one of: an edgeless bipartite graph, a complete bipartite graph, a complement of a perfect matching, a perfect matching, or a homogeneous generic bipartite graph. 
\end{lem}
\begin{proof}
If we do not ask $\Omega(G)$ to be finite, the proof of the locally finite analogue \cite[Lemma 4.5]{GMo} carries over. Compared to the locally finite case, we only have to deal with one other 'type' of graph, due to \cite[Remark 1.3]{GGK}
\end{proof}

\begin{lem}\label{lem_CHomRanBipImpliesHomRanBip}
Let $G$ be a C-homogeneous generic bipartite graph. Then $G$ is homogeneous bipartite.
\end{lem}

\begin{proof}
Let $VG=A\cup B$ be the natural bipartition of~$G$, let $X$ and $Y$ be two isomorphic induced finite subgraphs of~$G$, and let $\varphi:X\to Y$ be an isomorphism.
Let $a\in A\sm X$ be a vertex adjacent to all the vertices of $X\cap B$ and let $b\in B\sm X$ be a vertex adjacent to all the vertices of $X\cap A$ and to~$a$.
Let $a',b'$ be the corresponding vertices for $Y$. Since $G$ is bipartite, both $G[X+a+b]$ and $G[Y+a'+b']$ are connected induced subgraphs of~$G$ that are isomorphic to each other. Furthermore there is an isomorphism $\psi:G[X+a+b]\to G[Y+a'+b']$ such that the restriction of~$\psi$ to~$X$ is~$\varphi$. As there is an automorphism of~$G$ that extends $\psi$, this automorphism also extends $\varphi$ and $G$ is homogeneous.
\end{proof}

\begin{thm}\label{thm_BipartiteCHom}
A connected graph is a C-homogeneous bipartite graph if and only if it belongs to one of the following classes:
\begin{enumerate}[(i)]
\item $T_{\kappa, \lambda}$ for cardinals $\kappa, \lambda$;
\item $C_{2m}$ for $m \in \N$;
\item $K_{\kappa, \lambda}$ for cardinals $\kappa, \lambda$;
\item $CP_\kappa$ for a cardinal $\kappa$;
\item homogeneous generic bipartite graphs.
\end{enumerate}
\end{thm}

\begin{proof} The nontrivial part is to show that this list is complete. So consider an arbitrary connected C-homogeneous bipartite graph $G$ with bipartition $X \cup~Y$. If $G$ is a tree then it is obviously semi-regular and hence a $T_{\kappa,\lambda}$. So suppose $G$ contains a cycle. Then, since $G$ is C-homogeneous, each vertex lies on a cycle. Now $G$ is either a cycle, which is even since $G$ is bipartite, or at least one vertex in $G$ has a degree greater than $2$ and $G$ embeds a $C_4$, due to Lemma \ref{lem_EmbeddingC4}. Thus $\Omega(G)$ contains at least one edge and by Lemma \ref{lem_Omega} we have to consider the following cases:

\smallskip
\noindent {\bf Case 1:} {\it $\Omega(G)$ is complete bipartite}. Suppose that there is an induced path $P = uxyv$ in $G$. Then $\Omega(x,y)$ gives rise to an edge between $u$ and $v$, a contradiction. Hence $G$ is complete bipartite.

\smallskip
\noindent {\bf Case 2:} {\it $\Omega(G)$ is the complement of a perfect matching}. Consider $x \in X$ and $y \in Y$ such that $\{x,y\}$ is an edge of $G$. Since $\Omega(x,y)$ is the complement of a perfect matching and $G$ is not a cycle, there is an index set $I \supseteq \{1,2\}$ such that $N(x) = \{y\} \cup \left\{y_i | i \in I\right\}$, $N(y) = \{x\} \cup \left\{x_i | i \in I\right\}$ and for $i \in I$ the vertex $x_i$ is nonadjacent to $y_i$ but adjacent to all $y_j$ with $j \in I \sm \{i\}$. Since $\Omega(x, y_1)$ is also the complement of a perfect matching there is a unique vertex $a \in N(y_1) \sm N(y)$. Since $x_i$ with $i \neq 1$ is adjacent to $y_1$ it is contained in $\Omega(x,y_1)$ and therefore $y_i$ is adjacent to $a$. Thus for all $i \in I$ we have $N(y_i) = N(y) - x_i + a$. Now by symmetry there is a unique vertex $b$ adjacent to all $x_i$ with $i \in I$ but non-adjacent to $x$ and for all $i \in I$ there is $N(x_i) = N(x) - y_i + b$. If we look at $\Omega(x_1,y_2)$ we have $x,a \in N(y_2)$ and $y,b \in N(x_1)$ which implies $\{a,b\} \in EG$ and hence $N(a) = N(x) - y + b$ and $N(b) = N(y) - x + a$. Because $G$ is connected we have $X = N(y) + a$ and $Y = N(x) + b$ which means that $G$ is itself the complement of a perfect matching. 

\smallskip
\noindent {\bf Case 3:} {\it $\Omega(G)$ is a perfect matching}. For the same reason as for locally finite graphs this case cannot occur (cp.\ \cite[Theorem~4.6]{GMo}).

\smallskip
\noindent {\bf Case 4:} {\it $\Omega(G)$ is homogeneous generic bipartite}. Let $U$ and $W$ be two disjoint finite subsets of $X$ (of Y). Since $G$ is connected there is a finite connected induced subgraph $H \subset G$ that contains both $U$ and $W$. By genericity, we find an isomorphic copy $H_\Omega$ of $H$ in $\Omega(G)$. Because $G$ is C-homogeneous there is an automorphism $\varphi$ of $G$ with $H_\Omega^\varphi = H$. Now there is a vertex $v$ in $Y$ (in $X$) that is adjacent to all vertices in $U^{\varphi^{-1}}$ and non-adjacent to all vertices in $W^{\varphi^{-1}}$. Hence $v^\varphi$ is adjacent to all vertices in $U$ and none in $W$ which implies that $G$ is generic bipartite. Furthermore $G$ is homogeneous bipartite by Lemma~\ref{lem_CHomRanBipImpliesHomRanBip}, as it is C-homogeneous.
\end{proof}

\section{C-homogeneous digraphs of Type II}\label{sec_main}

It is well known that a transitive locally finite graph either contains one, two, or infinitely many ends.  
For arbitrary infinite graphs, this was proved by Diestel, Jung and Möller \cite{DJM}. 
Since the underlying undirected graph of a transitive digraph is also transitive, the same holds for infinite transitive digraphs.
The two-ended C-homogeneous digraphs have a very simple structure which we could easily derive from the results of the previous sections. But since two-ended connected transitive digraphs are locally finite \cite[Theorem~7]{DJM} we refer to Gray and Möller \cite[Theorem 6.2]{GMo} instead. Consequently, this section only deals with digraphs that have infinitely many ends.

As a first result we prove that no connected C-homogeneous digraph of Type II with more than one end contains any triangle.

\begin{lem}\label{lem_TypeIINoTriangle}
Let $D$ be a connected C-homogeneous digraph of Type II with more than one end.
Then $D$ contains no triangle.
\end{lem}

\begin{proof}
Let $\CF$ be a basic cut system.
Suppose that $D$ contains a triangle.
By Lemma~\ref{lem_WithTrianglesBlocksComplete} every $\CF$-blocks of~$D$ that contains an edge is a tournament.
Since each $\CF$-separator has to consist of precisely one vertex, each $\CF$-block contains edges, and the $\CF$-blocks have to be homogeneous tournaments.
Thus $D$ is of Type I in contradiction to the assumption.
\end{proof}

In preparation of the next lemma we introduce the following well-known construction: Given an edge-transitive bipartite digraph $\Delta$ with bipartition $A \cup B$ such that every edge is directed from $A$ to $B$ we define $DL(\Delta)$ to be the unique connected digraph such that each vertex separates the digraph, lies in exactly two  copies of $\Delta$, and has both in- and out-neighbors (cp.\ \cite{CPW,GMo}).

\begin{lem}\label{lem_SeparatorAndReachabilityDigraph}
Let $D$ be a connected C-homogeneous digraph of Type II with more than one end. If $D$ has connectivity $1$, then $D$ is isomorphic to $DL(\Delta(D))$.
\end{lem}

\begin{proof}
This is direct consequence of Lemma~\ref{lem_TypeIINoTriangle} and Lemma~\ref{lem_OneOut/InBlock}.
\end{proof}

In the next two theorems we prove that in the cases that the reachability digraph is either isomorphic to $CP_\kappa$ or to $K_{2,2}$ the digraph has connectivity at most $2$ and we determine the only digraphs with connectivity~$2$ and these reachability digraphs that might be C-homogeneous.

We first define a class of digraphs with connectivity $2$ and reachability digraph $CP_\kappa$. Given $2 \le m \in \N$ and a cardinal $\kappa \ge 3$ consider the tree $T_{\kappa,m}$ and let $U \cup W$ be its natural bipartition such that the vertices in $U$ have degree $m$. Now subdivide each edge once and endow the neighborhood of each $u \in U$ with a cyclic order. Then for each new vertex $y$ let $u_y$ be its unique neighbor in $U$ and denote by $\sigma(y)$ the successor of $y$ in $N(u_y)$. Then for each $w \in W$ and each $x \in  N(w)$ we add an edge directed from $x$ to all $\sigma(y)$ with $y \in N(w) - x$. Finally we delete the edges and vertices of the $T_{\kappa,m}$ to obtain the digraph $M(\kappa,m)$. The locally finite subclass of this class of digraphs coincides with those digraphs $M(k,n)$ for $k,n\in\nat$ that are described in~\cite[Section~5]{GMo}. In Figure~\ref{pic_M(3,3)} the digraph $M(3,3)$ is shown: once with its construction tree and once with its set of $\CF$-separators.

\begin{figure}[h]
\begin{center}
\includegraphics[width=.48\textwidth]{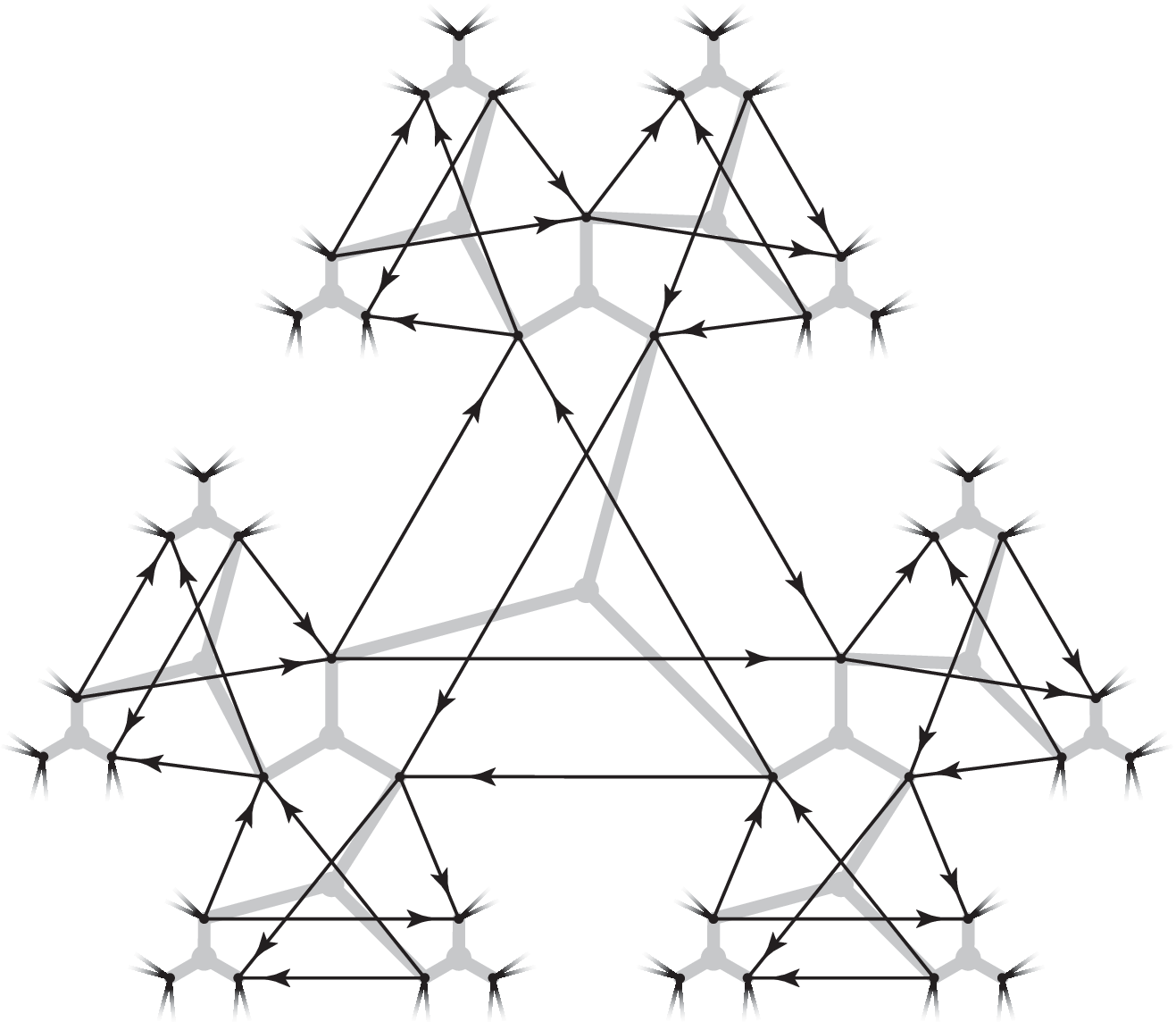}
\hfill
\includegraphics[width=.48\textwidth]{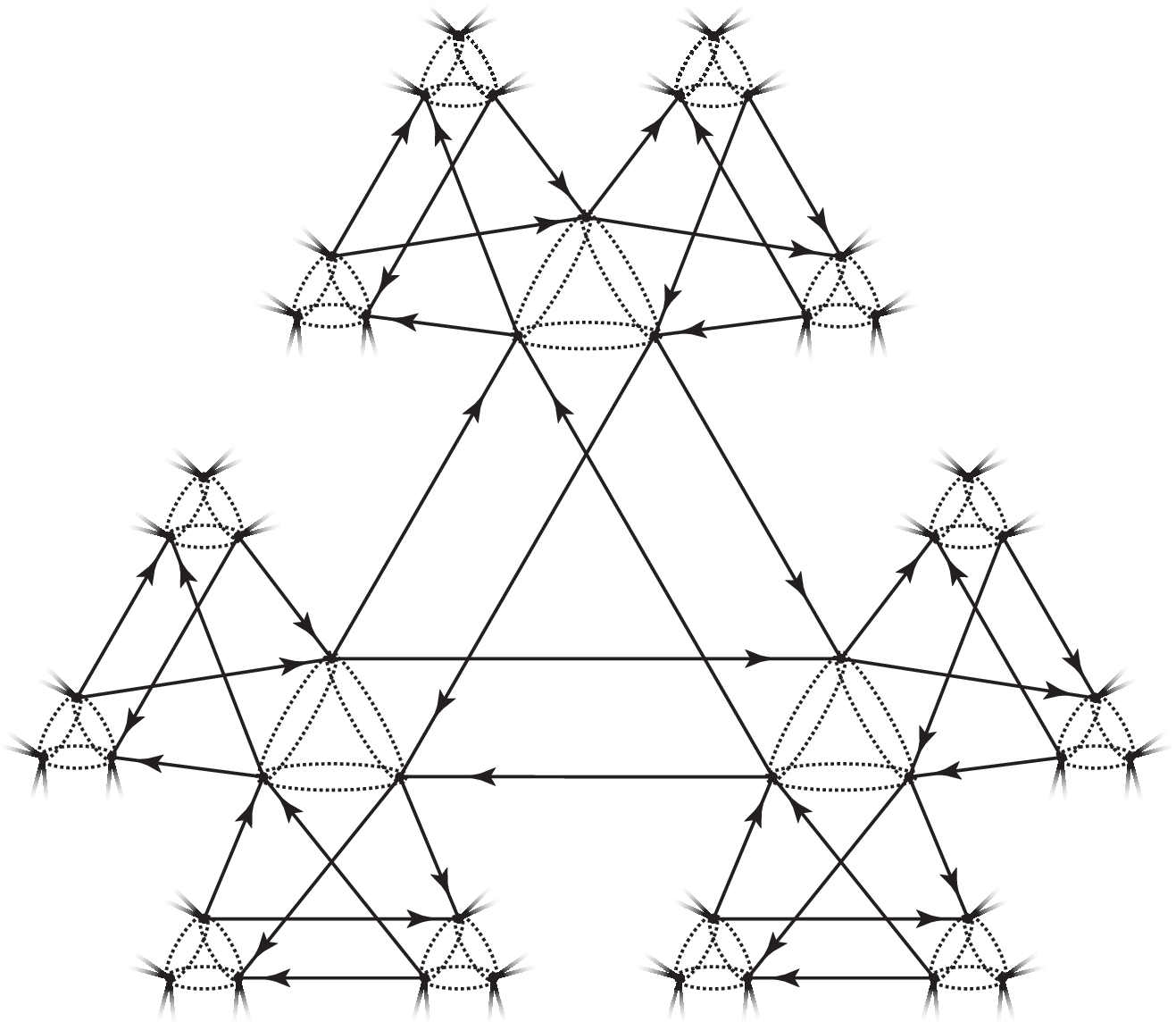}
\caption[Figure 1]{The digraph $M(3,3)$}\label{pic_M(3,3)}
\end{center}\end{figure}

\begin{thm}\label{thm_ComplementMatching}
Let $D$ be a connected C-homogeneous digraph of Type II with infinitely many ends and with $\Delta(D)\isom CP_\kappa$ for a cardinal $\kappa\ge 3$.
If $D$ has connectivity more than one, then $D$ is isomorphic to $M(\kappa,m)$ for an $m\in\nat$ with $m\ge 2$.
\end{thm}

\begin{proof}
By Lemma~\ref{lem_TypeIINoTriangle} the digraph $D$ contains no triangle.
Let $\CF$ be a basic cut system and let $\TF$ be the structure tree of~$D$ and $\CF$.
Let $S$ be a $\CF$-separator, let $X=\Delta_e$ for an $e \in ED$ such that $|S\cap X|\ge 2$, and let $K$ be a $\CF$-block with $S\sub K$ and $\Delta_e\sub K$, which all exists by Lemma~\ref{lem_SeparatorAndReachability}. Let $A\cup B$ be the natural bipartition of~$X$ such that its edges are directed from $A$ to~$B$.
For each $a\in A$ let us denote with $b_a$ the unique vertex in~$B$ such that $ab_a$ is no edge in~$X$. By symmetry we may assume that $A\cap S\ne\es$, so let $a\in A\cap S$.

First we will show that $X \cap S = \{a, b_a\}$. Since $S$ contains no edges by Lemma~\ref{lem_SeparatorEdgeless} it suffices to show that $A \cap S = \{a\}$. So let us suppose that there is another vertex $a' \neq a$ in $A\cap S$.
By C-homogeneity we have $A\sub S$.
Let $a' \in A$ be distinct to $a$ and $P$ an induced $a$-$a'$-path whose interior is contained in $D - K$. Denote the unique neighbor of $a$ on $P$ by $c$. Taking into account that $X$ is a $CP_\kappa$, there is a common successor for each pair of $A$-vertices; let $b$ be such a common successor of $a$ and $a'$. By C-homogeneity we can map $cPb$ onto $cPb_a$ by an isomorphism $\varphi$. Then $a^\varphi$ is a successor of $c$ that sends an edge to $b_a$. Hence $a^\varphi$ lies in $A$ and is distinct to $a$, contradicting the fact that $\desc{c}$ is a tree.

For the remainder let $X^0$, $S^0$ and $K^0$ refer to $X$, $S$ and $K$, respectively, and let $X^0 \cap S^0 = \{x_0, x_1\}$. Because each vertex clearly lies in exactly two distinct reachability digraphs, there is a unique reachability digraph $X^1 \neq X^0$ that contains $x_1$. If $x_0 \in X^1$ then it is straight forward to see that $D \isom M(\kappa, 2)$. So assume $x_0 \notin X^1$ and let $\psi$ be an automorphism of $D$ mapping $X^0$ onto $X^1$ and $x_0$ to $x_1$. Let $S^1$, $K^1$ denote the image under $\psi$ of $S^0$, $K^0$, respectively, and let $x_2 = x_1^\psi$. Since $\CF$ is basic there is an induced $x_0$-$x_1$-path $P$ the interior of which lies in $D - K^0$. We shall show that $P$ contains $x_2$.

Suppose that $P$ does not contain $x_2$ and has minimal length with this property. Let $u$ be the neighbor of $x_1$ on $P$, which clearly lies in $X^1$, and let $v$ be a neighbor of $u$ in $X^1$ distinct to $x_1$. If $v$ lies not on $P$, then $Puv$ is a path of the same length as $P$ which is induced by the minimality of $P$ and Theorem~\ref{thm_DescendentReachability}, contradicting the fact that $x_0$ and $v$ cannot lie in a common reachability digraph. On the other hand, if $v$ does lie on $P$ then consider a neighbor $w$ of $x_2$ in $X^1$ distinct to $v$. Remark that since $X^1$ is a $CP_\kappa$ there is an edge between $v$ and $x_2$. Thus by the choice of $P$ the path $Pvx_2w$ is induced and of the same length as $P$, which is impossible since $x_0$ and $w$ do not belong to a common reachability digraph. Hence $P$ contains $x_2$.

We have just proved that $\{x_1, x_2\}$ separates $x_0$ from any neighbor of $x_1$~in~$X^1$. Hence all $\CF$-separators have order $2$ and thus the blocks which contain edges consist each of a single reachability digraph. Now we repeat the previous construction to continue the sequences $(X^i)_{i \in \N}$, $(S^i)_{i \in \N}$, $(K^i)_{i \in \N}$ and $(x_i)_{i \in \N}$, respectively. Since $Px_2$ is an induced $x_0$-$x_2$-path the interior of which lies in $D - K^1$, we can apply the same argument as above to assure that $P$ contains $x_3$. Hence by induction we have $x_i \in P$ for all $i \in \N$, and since $P$ is finite there is an $m \in \N$ such that $x_m = x_0$. Furthermore we have $X^m = X^0$, $S^m = S^0$ and $K^m = K^0$. One can verify that $\{x_0, x_1, \dots, x_{m-1}\}$ forms a maximal $\CF$-inseparable set -- a $\CF$-block -- which means that $D$ is isomorphic to~$M(\kappa, m)$. 
\end{proof}

\begin{figure}[h]
\begin{center}
\includegraphics[width=.48\textwidth]{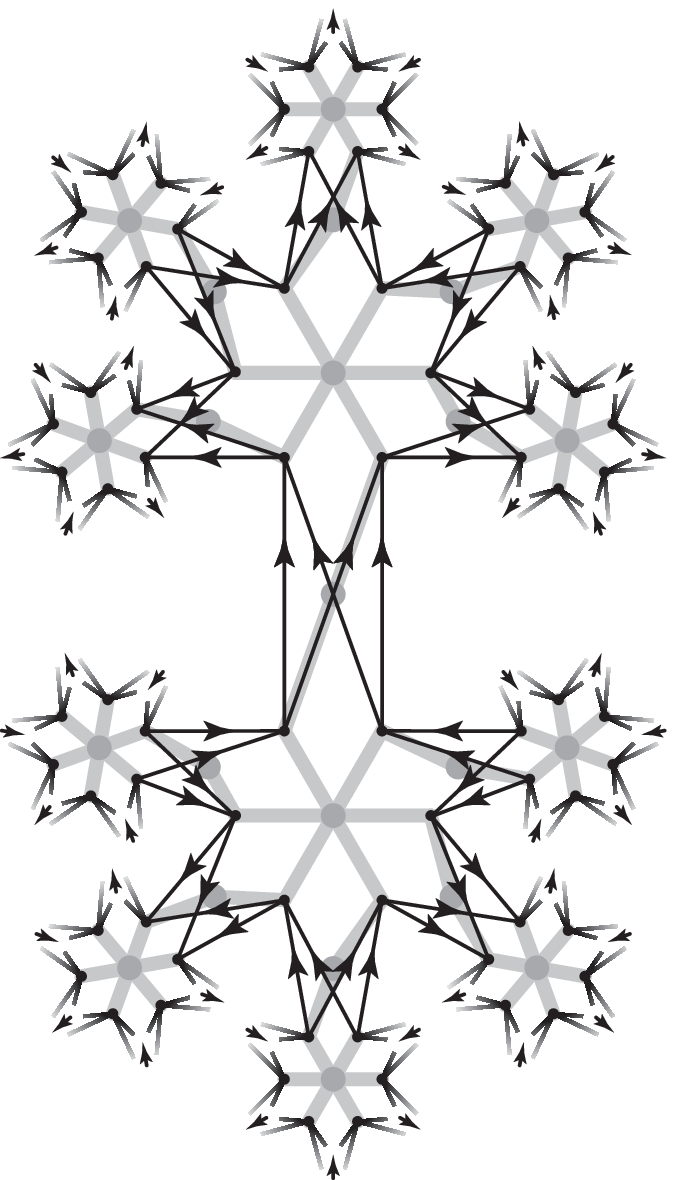}
\hfill
\includegraphics[width=.48\textwidth]{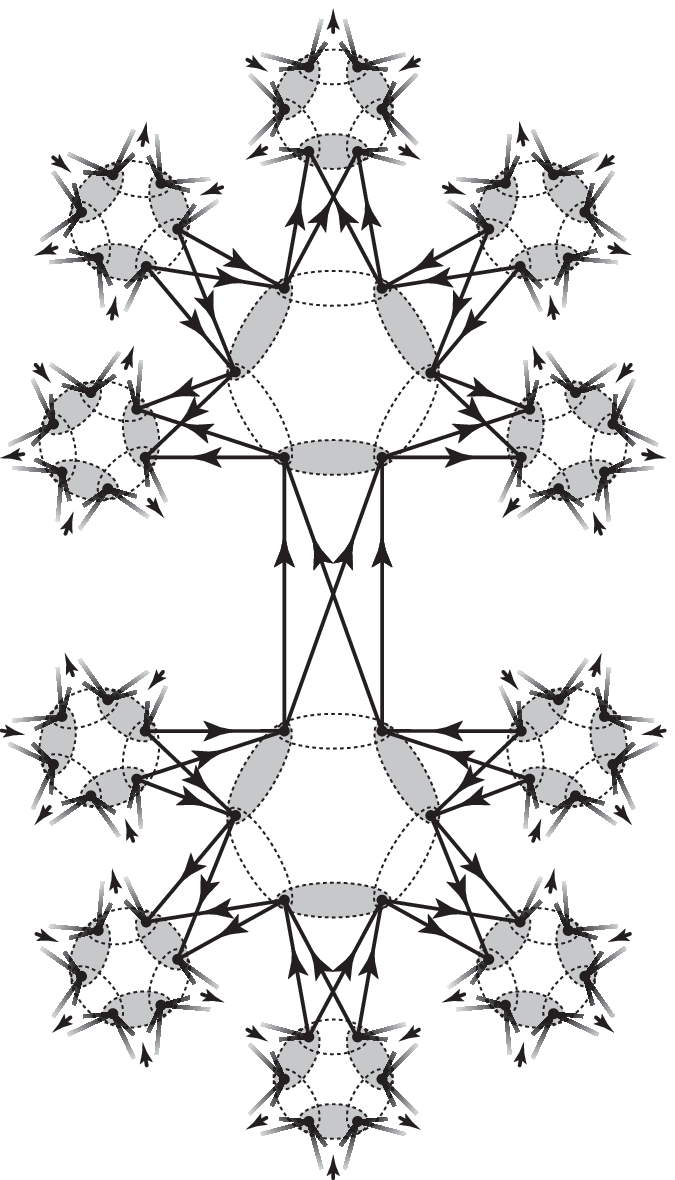}
\caption[Figure 1]{The digraph $M'(6)$}\label{pic_M'(6)}
\end{center}\end{figure}

In preparation of the next theorem we define a class of digraphs with connectivity $2$ and reachability digraph $K_{2,2}$. For $2 \le m \in \N$ consider the tree $T_{2,2m}$ and let $U \cup W$ be its natural bipartition such that the vertices in $U$ have degree $2m$. Now subdivide every edge once and enumerate the neighborhood of each $u \in U$ from $1$ to $2m$ in a such way that the two neighbors of each $w \in W$ have distinct parity. For each new vertex $x$ let $u_x$ be its unique neighbor in $U$ and define $\sigma(x)$ to be the successor of $x$ in the cyclic order of $N(u_x)$. For any $w \in W$ we have a neighbor $a_w$ with even index, and a neighbor $b_w$ with odd index. Then we add edges from both $a_w$ and $\sigma(a_w)$ to both $b_w$ and $\sigma(b_w)$. Finally we delete the $T_{2,2m}$. With $M'(2m)$ we denote the resulting digraph.
Figure~\ref{pic_M'(6)} shows the digraph $M'(6)$: on the left side with its construction tree and on the right side with the separators of the two possible basic cut systems.

\begin{thm}\label{thm_K22}
Let $D$ be a connected C-homogeneous digraph of Type II with infinitely many ends and with $\Delta(D)\isom K_{2,2}$.
If $D$ has connectivity more than one, then $D$ is isomorphic to $M'(2m)$ for $2 \le m\in\nat$.
\end{thm}

\begin{proof}
Lemma~\ref{lem_TypeIINoTriangle} implies that $D$ contains no triangle.
Let $\CF$ be a basic cut system of $D$.
Let $S^0$ be a $\CF$-separator and let $X^0=\Delta_e$ for an $e\in ED$ such that $|S^0\cap X^0|\ge 2$. Such an $X$ exists by Lemma~\ref{lem_SeparatorAndReachability}.
As $\Delta(D)\isom K_{2,2}$ and as no $\CF$-separator contains any edge by Lemma~\ref{lem_SeparatorEdgeless}, there is $|S^0\cap X^0|=2$.
So let $x_0,x_1$ be the two vertices in $X^0\cap S^0$. Let $X^1$ be the other reachability digraph that contains $x_1$ and let $x_2$ be the unique vertex in $X^1$ that is not adjacent to $x_1$. Let $\psi$ be an automorphism of $D$ that maps $X^0$ onto $X^1$ and let $S^1$ bet the image of $S^0$ under $\psi$. 

With the same technique as in the previous proof, we can verify that $\{x_1, x_2\}$ separates $D$, such that $S^0 = \{x_0, x_1\}$, we can continue the sequences $(x_i)_{i \in \N}$ and $(S^i)_{i \in \N}$, and there is $n \in \N$ such that $x_n = x_0$. Since $D$ has infinitely many ends we have $n \ge 3$, and as $x_i \in S^i$ only holds for all even integers $i$ we have $n = 2m$ with $m \ge 2$. Now analog as in the proof of Theorem~\ref{thm_ComplementMatching} $\bigcup_i{S^i}$ forms a $\CF$-block that contains edges. Hence there is only one $\Aut(D)$-orbit on the $\CF$-blocks and $D$ is isomorphic to $M'(2m)$.
\end{proof}

If we assume $\Delta(D)$ to be one of the other possibilities as described in Theorem~\ref{thm_BipartiteCHom}, then the C-homogeneous digraphs have - in contrast to the other two cases - connectivity~$1$.

\begin{lem}\label{lem_CycleCompleteRado}
Let $D$ be a connected C-homogeneous digraph of Type II with infinitely many ends and such that $\Delta(D)$ is isomorphic to a $T_{\kappa,\lambda}$ for cardinals $\kappa,\lambda$, a $C_{2m}$ with $4\leq m\in\nat$, a $K_{\kappa,\lambda}$ for cardinals $\kappa,\lambda\geq 2$, or an infinite homogeneous generic bipartite digraph. 
Then $D$ has connectivity~$1$.
\end{lem}

\begin{proof}
Since $D$ is of Type II, it contains no triangle by Lemma~\ref{lem_TypeIINoTriangle}.
Let us suppose that $D$ has connectivity at least~$2$ and let $\CF$ be a basic cut system of~$D$.
Let $S$ be a $\CF$-separator and let $X$ be a reachability digraph with $|S\cap X|\ge 2$ as in Lemma~\ref{lem_SeparatorAndReachability}.
We investigate the given reachability digraphs one by one and get in each case a contradiction and, thereby, we get a contradiction in general to the assumption that $D$ has connectivity at least~$2$.
So let us assume that $X\isom T_{\kappa,\lambda}$ for cardinals $\kappa,\lambda$.
By Lemma~\ref{lem_OneOut/InBlock} we know, that $\kappa,\lambda\ge 2$, as $D$ is no tree.
Let $x,y\in S\cap X$ such that $d_X(x,y)$ is maximal.
Such vertices exists as $S$ is finite.
Let $e_1$ be the first edge on the path from $x$ to~$y$ in~$X$ and let $e_2$ be another edge incident with~$x$.
Then there is an $\alpha\in\Aut(D)$ with $e_1^\alpha=e_2$.
But then $y^\alpha$ lies in a common separator with $x$, as $x^\alpha=x$.
By Corollary~\ref{cor_SeparatorToBlock} the separator $S^\alpha$ has to be the same as $S$.
But this contradicts the maximality of $d_X(x,y)$, as $d_X(y^\alpha,y)>d_X(x,y)$.

Let us now assume that $X\isom C_{2m}$ for a $4\leq m\in\nat$ and let $x,y$ be distinct vertices in $S\cap X$.
Then there is an induced path $P$ from $x$ to~$y$ that lies apart from $x$ and $y$ in a component of $D-S$ that intersects trivially with $X$.
We first show that we may assume that $d_X(x,y)\ge 4$.
Let $e_1,e_2$ be the two edges in $D[X]$ that are incident with $x$.
If $d_X(x,y)=k\le 3$, then let $\alpha\in\Aut(D)$ with $e_1^\alpha=e_2$.
Then there is $d_X(y,y^\alpha)=2k$, as $m\ge 4$.
Thus we have shown that there are $x,y\in S\cap X$ with $d_X(x,y)\ge 4$.
Let $s_1$ and $s_2$ be the vertices in~$X$ that are adjacent to~$y$ and let $t$ be a vertex in~$X$ that is adjacent to~$x$.
Since $d_X(x,y)\geq 4$, the graphs $txPys_i$ for $i=1,2$ are induced paths.
Hence there is an automorphism $\alpha$ of~$D$ that maps $txPys_1$ onto $txPys_2$ and thus $d_X(s_1,x)=d_X(s_2,x)$ and $d_X(s_1,t)=d_X(s_2,t)$, a contradiction as $X$ is a cycle.

For the next case let us assume that $X\isom K_{\kappa,\lambda}$ for cardinals $\kappa,\lambda\geq 2$.
Let $A\cup B$ be the natural bipartition of~$X$.
Since $|S\cap X|\ge 2$, the vertices in $S\cap X$ lie in the same set either $A$ or~$B$. So we may assume that they lie in~$A$.
By the C-homogeneity it is an immediate consequence that $A\sub S$.
As the $\CF$-separators have minimal cardinality with respect to separating ends, there is $|A|\le|B|$.
If there is a $\CF$-separator $S'$ with $|S'\cap B|\ge 2$, then $B\sub S'$.
If in addition the intersection of~$B$ with another reachability digraph distinct to~$X$ is~$B$, then it is a direct consequence that $\kappa=\lambda$ is finite and that $D$ has two ends.
Thus there are two distinct reachability digraphs $X_1,X_2$ that intersects with $B$ non-trivially and that are distinct to~$X$.
Let $A_1,B_1, A_2,B_2$ be the natural bipartitions of $X_1,X_2$, respectively.
Let $P$ be an induced path from $A_1\sm B$ to $A_2\cap B$ in a component of $D-S$ that intersects non-trivially with $X$.
Let $a$ be the vertex in $A\cap P$ that is adjacent to the vertex in $P\cap A_2$ and let $b$ be the vertex in $B\cap A_1$.
Then there is an automorphism $\alpha$ of~$D$ that maps $P$ onto $Pab$.
But this contradicts the fact that the endvertices of $Pab$ lie both in $A_1$ but the endvertices of $P$ do not lie in the same first component of any reachability digraph as $A_1\ne A_2$.
Thus we conclude that $|B\cap S'|=1$.
So let $x,y,z\in B$ be three distinct vertices.
There is a shortest induced path $P$ from $x$ to $y$ in that component of $D-S$ that contains $B$.
Let $a\in A$ and let $b$ be the vertex on~$P$ with distance~$2$ to~$y$.
Then there is an automorphism $\alpha$ of~$D$ that maps $zaxPb$ onto $yaxPb$.
Thus we conclude that $d(b,z)=2$. But then $z$ has to have incident edges that are directed both towards or both from distinct $\CF$-blocks. This contradicts Lemma~\ref{lem_OneOut/InBlock}.

Let us finally assume that $X$ is isomorphic to an infinite homogeneous generic bipartite digraph.
Let again $A\cup B$ be the natural bipartition of~$X$.
Since $X$ is homogeneous, all vertices in the same set $A$ or $B$ have distance $2$ to each other.
We conclude that $|S\cap A|\ge 2$ immediately implies $A\sub S$ which contradicts the finiteness of~$S$.
Conversely we also know $|B\cap S|\le 1$.
Since $D$ has connectivity at least $2$, there is $|A\cap S|=1=|B\cap S|$.
Let $a,b$ be the vertices in $A\cap S,B\cap S$, respectively, and let $ab'a'b$ be a path of length $3$ from $a$ to $b$. This path exists because each two vertices in the same set $A$ or $B$ have distance $2$ to each other as before.
Since there are infinitely many vertices in~$A$ that are adjacent to~$b'$ but not to~$b$, all these vertices have to lie in~$S$, a contradiction.
Thus we conclude that $D$ has connectivity~$1$.
\end{proof}

Let us summarize the conclusions of this section in the following theorem.
In its proof we will finally prove that all the candidates for C-homogeneous digraphs are really C-homogeneous.

\begin{thm}\label{thm_Main}
Let $D$ be a connected digraph of Type II with infinitely many ends.
Then $D$ is C-homogeneous if and only if one of the following holds:
\begin{enumerate}[(1)]
\item $\Delta(D)\isom CP_\kappa$ for a cardinal $\kappa\geq 3$ and $D\isom DL(\Delta(D))$.
\item $\Delta(D)\isom C_{2m}$ for $2\leq m\in\nat$ and $D\isom DL(\Delta(D))$.
\item $\Delta(D)\isom K_{\kappa,\lambda}$ for cardinals $\kappa,\lambda\geq 2$ and $D\isom DL(\Delta(D))$.
\item $\Delta(D)$ is isomorphic to an infinite homogeneous generic bipartite digraph and $D\isom DL(\Delta(D))$.
\item $\Delta(D)=CP_\kappa$ and $D\isom M(\kappa,m)$ for a cardinal $\kappa \ge 3$ and $2 \le m\in\nat$.
\item $\Delta(D)=K_{2,2}$ and $D\isom M'(2m)$ for $2 \le m\in\nat$.
\end{enumerate}
\end{thm}

\begin{proof}
By the Lemmas~\ref{lem_TypeIINoTriangle}, \ref{lem_SeparatorAndReachabilityDigraph}, and \ref{lem_CycleCompleteRado} and by the Theorems~\ref{thm_ComplementMatching} and \ref{thm_K22}, it remains to show that the described digraphs are indeed C-homogeneous. Remark that the underlying undirected graph of $DL(T_{\kappa,\lambda})$ is a regular tree and thus $DL(T_{\kappa,\lambda})$ is not of Type II.
It is straight forward to see that the graphs of the part (1)-(4) are C-homogeneous. So let $D\isom M(\kappa,m)$ for an $m\in\nat$ with $m\ge 2$ and a cardinal $\kappa$.
Let $\CF$ be a basic cut system of~$D$.
Let $A$ and $B$ be two connected induced finite and isomorphic subdigraphs of~$D$ and let $\varphi$ be an isomorphism from $A$ to~$B$.
Let us first consider the case that $A$ contains no $2$-arc. Then both $A$ and $B$ lie in a reachability digraph, each.
Without loss of generality we may assume that they lie in the same reachability digraph $\Delta$ of~$D$.
But, as the reachability-digraphs are obviously C-homogeneous, it is straight forward to see that the isomorphism $\varphi$ from $A$ to~$B$ first extends to an automorphism of~$\Delta$ and then also to an automorphism of~$D$.
So let us assume that $A$ contains a $2$-arc.
Let $S$ be a $\CF$-separator such that $A\sm S$ has at least two components.
Since $\abs{S}\le 2$, for at least one of the components of $D\sm S$, let us denote this with $K$, there is a connected subdigraph $A_1$ of $A$ such that $A\cap S=A_1\cap S$ and $A\cap K=A_1\sm S$.
Then there is a $\CF$-separator $S_B$ such that $(A\cap S)^\varphi=B\cap S_B$.
By induction we can extend $\varphi|_{A_1}$ to an automorphism $\psi^*$ of~$D$. We will define the automorphism $\psi$ of~$D$ step by step. So let $\psi|_{K\cup S}:=\psi^*|_{K\cup S}$. Then $(K\cup S)^\psi$ is precisely $S_B$ together with that component of $D\sm S_B$ that contains $(A_1\sm S)^\varphi$. So we just have to define $\psi$ on that component of $D-S$ other than $K$. If $A_2:=A\sm(A_1\sm S)$ is connected, it is an inductive argument that we can define $\psi$ analog on the other component of $D\sm S$ as we did it on $C$.
Thus we may assume that $A_2$ consists of at least two components but then it has precisely two components as $\abs{S}=2$.
Let $X$ be the first $\CF$-block adjacent to~$S$ in the structure tree $\TF$ of~$\CF$ and $D$, that lies in that component of~$\TF-S$ that intersects non-trivially with $A_2\sm S$.
We distinguish between the two cases that $X$ contains edges or does not contain any edge of~$D$.

Let us first consider the case that $X$ contains edges of~$D$. Then $X$ is a reachability digraph of~$D$ (and hence isomorphic to a $CP_\kappa$) and $A\cap X$ must consist of precisely two edges.
The same must be true for the $\CF$-block $Y$ adjacent to~$S_B$ that has the same role as~$X$ just for $B$ instead of~$A$. Thus we can extend our definition of~$\psi$ from $X$ to~$Y$ and also to all components of $D-X$ that intersects with $A_2$ trivially. By induction on $\abs{A_2}$ we have constructed $\psi$.
Thus we may assume that $X$ does not contain any edge.
There is an enumeration $x_1,\ldots,x_m$ of the vertices of~$X$ such that $\{x_m,x_1\}$ and for all $i\leq m$ also $\{x_i,x_{i+1}\}$ are all the $\CF$-separators in~$X$. We may assume that $S=\{x_1,x_2\}$.
As the image of~$S$ under $\psi$ has already been defined, it is an immediate consequence, that we can define $\psi$ also inductively the other vertices $x_i$ and thus on~$X$. By induction we have defined $\psi$ on all components of $D-X$ and thus there is an automorphism of~$D$ that extends $\varphi$.

In the case that $D\isom M'(2m)$ for an $m\in\nat$ the arguments used are analog ones as in the case $D\isom M(\kappa,m)$ and therefore we omit that proof here.
\end{proof}

It is well known (see \cite{CPW}) that line digraphs of highly-arc-transitive digraphs are again highly-arc-transitive.
In some cases also C-homogeneity is preserved under taking the line digraph: Gray and M\"oller \cite{GMo} stated that the line digraph of a $DL(C_{2m})$ is C-homogeneous.
In terms of our classification:
\begin{rem}\label{rem_LineOfCycle}
For each $m\in \N$ we have $L(DL(C_{2m}))\isom M'(2m)$.
\end{rem}

\begin{proof}
Consider the digraph $D = DL(C_{2m})$ for a $m \in \N$. By construction the deletion of each single vertex $v$ of $D$ splits the digraph into two components such that $v$ has two out-neighbors in the one and two in-neighbors in the other component. Thus the four edges that are incident with $v$ form a $K_{2,2}$ in $L(D)$ whose independent vertex sets separate $L(D)$. Furthermore the edges of each $C_{2m}$ in $D$ form an independent set in $L(D)$ so that any two adjacent edges lie in a common $K_{2,2}$ in $L(D)$. One can easily verify that this digraph is indeed isomorphic to $M'(2m)$. 
\end{proof}

Interestingly, our classification implies that C-homogeneity is not generally preserved under taking line digraphs.
Indeed, For all $m \in \N$ the line digraph of $M'(2m)$ is triangle-free, has infinitely many ends, and has connectivity~$4$, hence it is not of Type II. 
Thus, by Theorem \ref{thm_Main}, we know that $L(M'(2m)) \isom L(L(DL(C_{2m})))$ is not C-homogeneous.
This had remained an open question in~\cite{GMo}.

\section{Final remarks}\label{sec_Final}

Let us take a closer look at two specific kinds of digraphs that occur as `building blocks' in our classification.
The first kind are the homogeneous tournaments, which feature in our classification of the connected C-homogeneous digraphs of Type~I.
While Lachlan \cite{L} classified the countable homogeneous tournaments, no characterization is known for the uncountable ones.
The second kind of building blocks that deserve a closer look are the generic homogeneous bipartite graphs, which occur in the classification of the connected C-homogeneous digraphs of Type II.
There is exactly one countable such digraph (\cite[Fact~1.2]{GGK}), but it is shown in~\cite{GGK} that the number of isomorphism types of homogeneous generic bipartite graphs with $\aleph_0$ vertices on the one side of the bipartition and $2^{\aleph_0}$ vertices on the other side is independent of ZFC.
Hence, classifying the uncountable generic homogeneous bipartite graphs remains an undecidable problem.

\end{document}